\documentclass{article}

\usepackage{arxiv}

\usepackage[utf8]{inputenc} % allow utf-8 input
\usepackage[T1]{fontenc}    % use 8-bit T1 fonts
\usepackage{hyperref}       % hyperlinks
\usepackage{url}            % simple URL typesetting
\usepackage{booktabs}       % professional-quality tables
\usepackage{amsfonts}       % blackboard math symbols
\usepackage{nicefrac}       % compact symbols for 1/2, etc.
\usepackage{microtype}      % microtypography
\usepackage{lipsum}
\usepackage{graphicx}
\usepackage{amsmath}
\graphicspath{ {./images/} }
\usepackage{caption}
\usepackage{subcaption}
\usepackage{lscape}

\title{A  multifidelity approach to continual learning for physical systems}

\author{
Amanda Howard, Yucheng Fu, and Panos Stinis \\ Advanced Computing, Mathematics and Data Division \\
Pacific Northwest National Laboratory\\
Richland, WA 99354 \\
  \texttt{amanda.howard@pnnl.gov} 
}

\begin{document}
\maketitle
\begin{abstract}
We introduce a novel continual learning method based on multifidelity deep neural networks. This method learns the correlation between the output of previously trained models and the desired output of the model on the current training dataset, limiting catastrophic forgetting. On its own the multifidelity continual learning method shows robust results that limit forgetting across several datasets. Additionally, we show that the multifidelity method can be combined with existing continual learning methods, including replay and memory aware synapses, to further limit catastrophic forgetting. The proposed continual learning method is especially suited for physical problems where the data satisfy the same physical laws on each domain, or for physics-informed neural networks, because in these cases we expect there to be a strong correlation between the output of the previous model and the model on the current training domain.  
\end{abstract}

\begin{figure}[h]
    \centering
    \includegraphics[width=\textwidth]{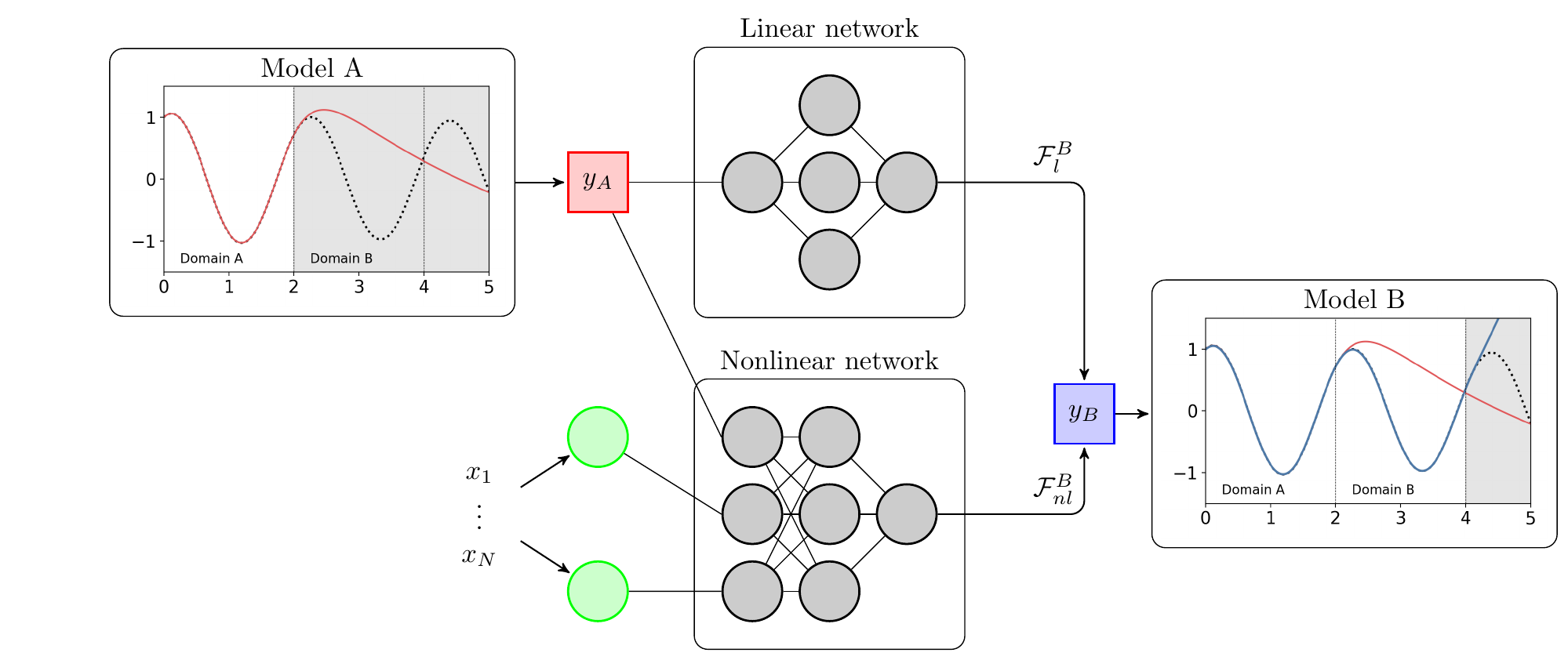}
    \caption{Graphical abstract}
    \label{fig:Graphical_abstract}
\end{figure}
% keywords can be removed
%\keywords{First keyword \and Second keyword \and More}

\section{Introduction}

In many real world applications of machine learning data is received sequentially or in discrete datasets. When used as training data new information received about the system requires completely retraining a given neural network. Much recent work has focused on how to instead incorporate the newly received training data into the machine learning model without requiring retraining with the full dataset and without forgetting the previously learned model. This process is referred to as continual learning \cite{parisi2019continual}. One key goal in continual learning is to limit catastrophic forgetting, or abruptly and completely forgetting the previously trained data. 

Many methods have been proposed to limit forgetting in continual learning. In replay (rehearsal), a subset of the training set from previously trained regions are used in training subsequent models, so the method can limit forgetting by reevaluating on the previous regions \cite{verwimp2021rehearsal}. However, replay requires access to the previously used training data sets. This both requires large storage capabilities for large datasets, and also physical access to the previous dataset. However, data privacy can limit access to prior datasets, so replay may not be a feasible option. 
 An alternative to replay are regularization methods, where a regularizer is used to assign weights to each parameter in the neural network, representing the parameter's importance. Then, a penalty is applied to prevent the parameters with the largest weights from changing. Multiple methods have been proposed for how to calculate the importance weights. Among the top choices are Synaptic Intelligence \cite{zenke2017continual}, elastic weight consolidation (EWC) \cite{kirkpatrick2017overcoming}, and memory aware synapses (MAS) \cite{aljundi2018memory}. Subsequent work has shown that MAS performs among the best in multiple use cases, and is more robust to the choice of hyperparameters, so here we use MAS \cite{de2021continual, hsu2018re}. Finally, a third category of continual learning methods includes those that employ task-specific modules \cite{rusu2016progressive}, ensembles \cite{wen2020batchensemble}, adapters \cite{pfeiffer2020adapterhub}, reservoir computing based architectures \cite{bereska2022continual}, slow-fast weights \cite{munkhdalai2017meta,vladymyrov2023continual} and more.

%\cite{tercan2022continual} uses MAS, initializes with lowest previous error on new training set. So for $\Omega_D$, evaluate A, B, C on $\Omega_D$ and use weights with lowest error as initialization. This is a straight forward thing to try--but weird with the multifidelity. 

In recent years, a huge research focus has been on scientific machine learning methods for physical systems \cite{karniadakis2021physics, baker2019workshop, cuomo2022scientific}, for example fluid mechanics and rheology \cite{jin2021nsfnets, raissi2020hidden, cai2021physics, Joglekar_2023}, metamaterial development \cite{liu2019multi, chen2020physics, fang2019deep}, high speed flows \cite{mao2020physics}, and power systems \cite{misyris2020physics, huang2022applications, moya2023dae}.  In particular, physics-informed neural networks, or PINNs \cite{raissi2019physics}, allow for accurately representing differential operators through automatic differentiation, allowing for finding the solution to PDEs without explicit mesh generation. Work on continual learning for PINNs is limited. While as a first attempt PINNs can be trained on the entire domain because the issues of data acquisition and privacy do not apply, many systems have been identified for which it is not possible to train a PINN for the entire desired time domain. For example, even the simple examples used in this work, a pendulum and the Allen-Cahn equation, cannot be trained by a PINN for long times. Recent work has looked at improving the training of PINNs for such systems, including applications of the neural tangent kernel \cite{wang2023long}, but more work remains to be done. The closest work we are aware of for continual learning with PINNS is the backward-compatible PINNs in \cite{mattey2022novel} and incremental PINNs (iPINNs) in \cite{dekhovich2023ipinns}.  Backward-compatible PINNs train $N$ PINNs on a sequence of $N$ time domains, and in each new domain enforce that the output from the current PINN satisfies the PINN loss function in the current domain and the output of the previous model on all previous domains. We note that this work is distinct from the replay approach taken with PINNs in this work, both in the single fidelity and multifidelity cases, because we enforce that the $N$th neural network satisfies the residual in all prior domains, not the output from the previous model. In iPINNs, PINNs are trained to satisfy a series of different equations through a subnetwork for each equation, rather than the same equation over a long time. 

%In contrast to \cite{shukla2021parallel}, the solutions in each domain are learned sequentially instead of in parallel. \cite{meng2020ppinn} "hence decomposing a long-time problem into many independent short-time problems supervised by an inexpensive/fast coarse-grained (CG) solver. In particular, the serial CG solver is designed to provide approximate predictions of the solution at discrete times, while initiate many fine PINNs simultaneously to correct the solution iteratively."

We will introduce the multifidelity continual learning method in Sec. \ref{sec:method}. We will then show the performance of the method on physics-informed problems in Sec. \ref{sec:pinns} and on data-informed problems in Sec. \ref{sec:data}.

\section{Multifidelity continual learning method}\label{sec:method}

We assume that we have a domain $\Omega$, which we divide into $N$ subdomains $\Omega = \cup_{i=0}^N \Omega_i$. We will learn sequential models on each subdomain $\Omega_i$, with the goal that the $i$th model can provide accurate predictions on the domain $\cup_{j=0}^i \Omega_j$. That is, the $i$th model does not forget the information learned on earlier domains used in training. We will focus on applications to physical systems, where we either have data available or knowledge of the physical laws the system obeys. We will begin this section with a brief overview of physics-informed neural networks (PINNs), then discuss the multifidelity continual learning method (MFCL), and conclude with a description of methods we use to limit catastrophic forgetting. 

\subsection{Physics-informed neural networks}

In this section we give a brief introduction to single-fidelity and multifidelity physics-informed neural networks (PINNs), which were introduced in \cite{raissi2019physics} and have been covered in depth for many relevant applications \cite{rasht2022physics, karniadakis2021physics}. PINNs are generally used, in these applications, for initial-boundary valued problems.
\begin{align}
    &\mathbf{s}_t + \mathcal{O}_{\mathbf{x}} [ \mathbf{s}] = \mathbf{0}, \; \mathbf{x} \in \Omega, t\in[0, T] \\
    &\mathbf{s}(\mathbf{x}, t) = \mathbf{g}(\mathbf{x}, t)  \; \mathbf{x} \in \partial \Omega, t\in[0, T] \\
    &\mathbf{s}(\mathbf{x}, 0) = \mathbf{u}(\mathbf{x})  \; \mathbf{x} \in  \Omega
\end{align}
where $\Omega \in \mathbb{R}^N$ is an open, bounded domain with boundary $\partial \Omega$, $\mathbf{g}$ and $\mathbf{u}$ are given functions, and $\mathbf{x}$ and $t$ are the spatial and temporal coordinates, respectively. $\mathcal{O}_{\mathbf{x}}$ is a general differential operator with respect to $\mathbf{x}$. We wish to find an approximation to $\mathbf{s}(\mathbf{x}, t)$ by a (series) of deep neural networks with parameters $\gamma$, denoted by $\mathbf{s}^\gamma(\mathbf{x}, t)$. The neural network is trained by minimizing the loss function 
\begin{equation}
    \mathcal{L}(\gamma) = \lambda_{bc} \mathcal{L}_{bc}(\gamma) + \lambda_{ic} \mathcal{L}_{ic}(\gamma) + \lambda_{r} \mathcal{L}_{r}(\gamma) + \lambda_{data} \mathcal{L}_{data}(\gamma)
\end{equation}
where the subscripts \emph{bc}, \emph{ic}, \emph{r}, and \emph{data} denote the terms corresponding to the boundary conditions, initial conditions, and residual, and any provided data, respectively. We take $N_{bc}$, $N_{ic}$, and $N_{r}$ to be the batch sizes of the boundary, initial, and residual data point, and denote the training data by $\left\{(\mathbf{x}_{bc}^i, t_{bc}^i), \mathbf{g}(\mathbf{x}_{bc}^i, t_{bc}^i) \right\}_{i=0}^{N_{bc}}$, $\left\{(\mathbf{x}_{ic}^i), \mathbf{u}(\mathbf{x}_{bc}^i) \right\}_{i=0}^{N_{ic}}$, and $\left\{(\mathbf{x}_{r}^i, t_{r}^i)\right\}_{i=0}^{N_{r}}$. The boundary and initial collocation points are randomly sampled uniformly in their respective domains. The selection of the $N_{r}$ residual points will be discussed in Sec. \ref{sec:rdps}. If data representing the solution $\mathbf{s}$ is available, we can also consider an additional dataset $\left\{(\mathbf{x}_{data}^i, t_{data}^i), \mathbf{s}(\mathbf{x}_{data}^i, t_{data}^i) \right\}_{i=0}^{N_{data}}$. This term is included to capture the data-based training we will cover in Sec. \ref{sec:data}. 

The individual loss terms are given by the mean square errors, 
\begin{align}
&\mathcal{L}_{bc}(\gamma) = \frac{1}{N_{bc}}  \sum_{i=0}^{N_{bc}} \left| \mathbf{s}_\gamma(\mathbf{x}_{bc}^i, t_{bc}^i)- \mathbf{g}(\mathbf{x}_{bc}^i, t_{bc}^i) \right|^2 \\
&\mathcal{L}_{ic}(\gamma) =  \frac{1}{N_{ic}} \sum_{i=0}^{N_{ic}} \left| \mathbf{s}_\gamma(\mathbf{x}_{ic}^i, 0)- \mathbf{u}(\mathbf{x}_{ic}^i) \right|^2 \\
&\mathcal{L}_{r}(\gamma) =  \frac{1}{N_{r}} \sum_{i=0}^{N_{r}} \left| \mathbf{r}_\gamma(\mathbf{x}_{r}^i, t_{r}^i)\right|^2
 \\
&\mathcal{L}_{data}(\gamma) =  \frac{1}{N_{data}} \sum_{i=0}^{N_{data}} \left| \mathbf{s}_\gamma(\mathbf{x}_{data}^i, t_{data}^i)- \mathbf{s}(\mathbf{x}_{data}^i, t_{data}^i) \right|^2 
\end{align}
where
\begin{equation}
    \mathbf{r}_\gamma(\mathbf{x}, t) = \frac{\partial}{\partial t} \mathbf{s}_\gamma(\mathbf{x}, t) + \mathcal{O}_\mathbf{x} [ \mathbf{s}_\gamma(\mathbf{x}_{data}^i, t_{data}^i)].
\end{equation}
The weighting parameters $\lambda_{bc}$, $\lambda_{ic}$, $\lambda_{r}$, and $\lambda_{data}$ are chosen before training by the user. 

Multifidelity PINNs, as used in this work, are inspired by \cite{meng2020composite}. We assume we have a low fidelity model in the form of a deep neural network that approximates a given dataset or differential operator with low accuracy. We want to train two additional neural networks to learn the linear and nonlinear correlations between the low fidelity approximation and a high fidelity approximation or high fidelity data. We denote these neural networks as $\mathcal{NN}_{l}$ for the linear correlation and $\mathcal{NN}_{nl}$ for the nonlinear correlation. The output is then $\mathbf{s}_\gamma(\mathbf{x}, t) = \mathcal{NN}_{nl}(\mathbf{x}, t; \gamma) + \mathcal{NN}_{l}(\mathbf{x}, t; \gamma) $, where $\gamma$ is all trainable parameters of the linear and nonlinear networks. The loss function includes an additional term, 
\begin{equation}
    \mathcal{L}_{MF}(\gamma) = \lambda_{bc} \mathcal{L}_{bc}(\gamma) + \lambda_{ic} \mathcal{L}_{ic}(\gamma) + \lambda_{r} \mathcal{L}_{r}(\gamma) + \lambda_{data} \mathcal{L}_{data}(\gamma) + \lambda \sum (\gamma_{nl, ij})^2, \label{eq:loss_MF}
\end{equation}
where $\{ \gamma_{nl, ij}\}$ is the set of all weights and biases of the nonlinear network $\mathcal{NN}_{nl}$. No activation function is used in $\mathcal{NN}_{l}$ to result in learning a linear correlation between the previous prediction and the high fidelity model.

\subsection{Multifidelity continual learning}
\begin{figure}[h]
    \centering
    \includegraphics[width=0.7\textwidth]{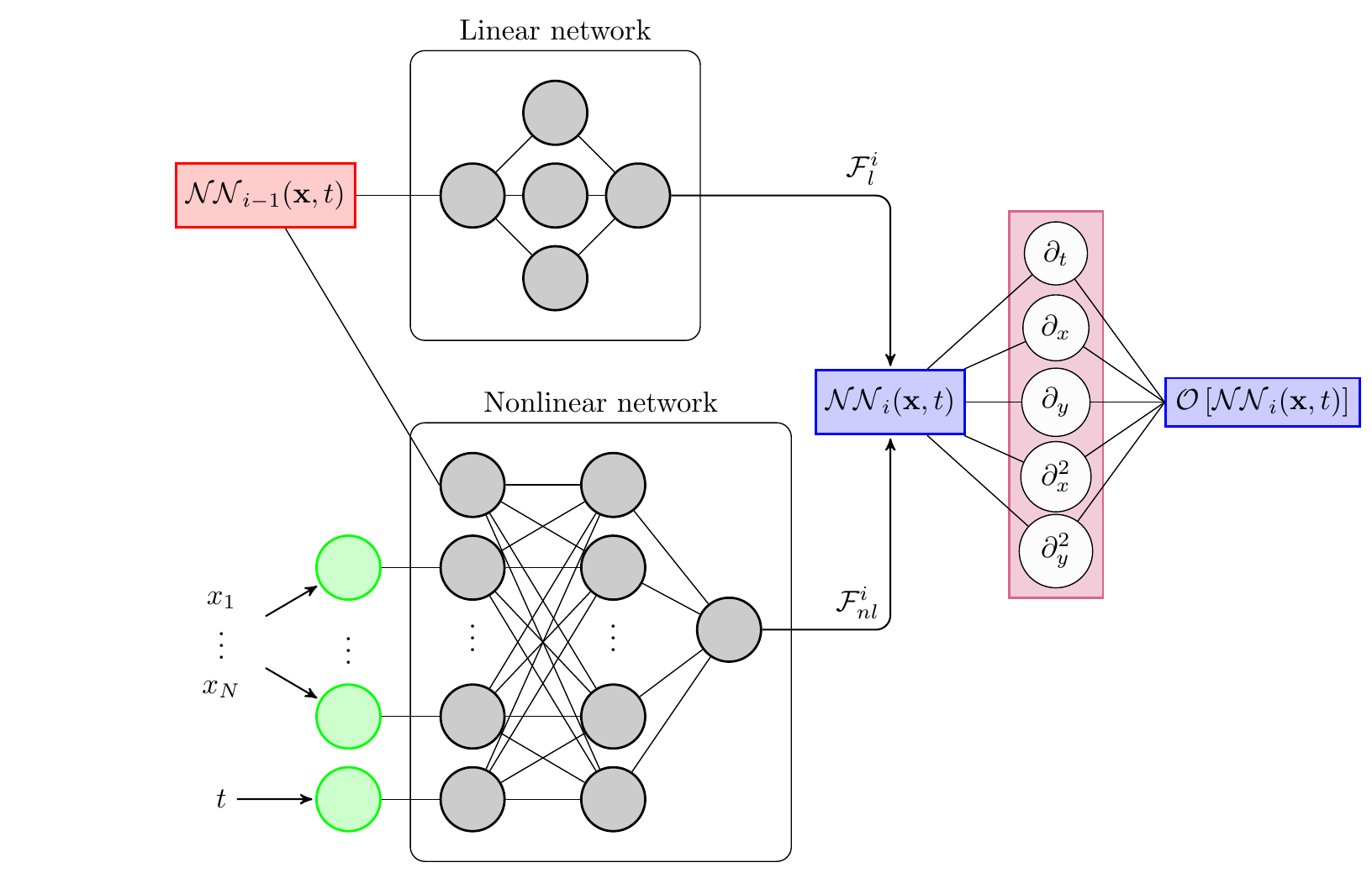}
    \caption{Diagram of the MF-CL method on domain $\Omega_i$. The output from the previously trained neural network, $\mathcal{NN}_{i-1}(\mathbf{x}, t; \gamma_{i-1})$, is used as input to the linear and nonlinear subnets for a point $(\mathbf{x}, t) \in \Omega_i$, $\mathbf{x}\in \mathbb{R}^N$. The output neural network is the sum of the linear and nonlinear subnetworks.}
    \label{fig:MF_NN_diagram}
\end{figure}

In the MFCL method, we exploit correlations between the previously trained models on prior domains and the expected model on the current domain. \emph{Explicitly, we use the prior model $\mathcal{NN}_{i-1}$ as a low fidelity model for domain $\Omega_i$. Then, we learn the correlation between $\mathcal{NN}_{i-1}$ on domain $\Omega_i$ and the data or physics given on the domain}. By learning a general combination of linear and nonlinear terms, we can capture complex correlations. Because the method learns only the correlation between the previous model and the new model, we can in general use smaller networks in each subdomain. The procedure requires two initial steps:
\begin{enumerate}
    \item Train a (single-fidelity) DNN or PINN on $\Omega_1$, denoted by $\mathcal{NN}^*(\mathbf{x}, t; \gamma^*)$. This network will approximate the solution in a single domain. 
    \item Train a multifidelity DNN or PINN in $\Omega_1$, which takes as input the single fidelity model $\mathcal{NN}^*(\mathbf{x}, t; \gamma^*)$ as a low fidelity approximation. This initial multifidelity network is denoted by $\mathcal{NN}_1(\mathbf{x}, t; \gamma_1)$
\end{enumerate}
Then, for each additional domain $\Omega_i$, we train  a multifidelity DNN or PINN in $\Omega_i$, denoted by $\mathcal{NN}_{i}(\mathbf{x}, t; \gamma_{i})$, which takes as input the previous multifidelity model $\mathcal{NN}_{i-1}(\mathbf{x}, t; \gamma_{i-1})$ as a low fidelity approximation. The goal is for $\mathcal{NN}_{i}(\mathbf{x}, t; \gamma_{i})$ to provide an accurate solution on $\cup_{j=1}^i \Omega_i$, even when data from $\Omega_j$, $j< i$, is not used in training the multifidelity network $\mathcal{NN}_{i}$. A diagram of the method is given in Fig. \ref{fig:MF_NN_diagram}. 

As we will show, the MF-CL method provides more accurate results with less forgetting than single fidelity training on its own, however, the method can be improved by a few methods that have been previously developed both for reducing forgetting in continual learning and for selecting collocation points for training PINNs. These methods are discussed below.

\subsection{Memory aware synapses}

Memory aware synapses (MAS) is a continual learning method that attempts to limit forgetting in continual learning by assigning an importance weight to each neuron in the neural network. Then, a penalty term is added to the loss function to prevent large deviations in the values of important weights when the next networks are trained. The importance weights are found by measuring how sensitive the output of neural net $\mathcal{NN}_n$ is to changes in the network parameters \cite{aljundi2018memory}. For each weight and bias $\gamma_{ij}$ in the neural network we calculate the importance weight parameter
\begin{equation}
    \Omega^n_{ij} = \frac{1}{N}\sum_{k=1}^N \left\| \frac{\partial\left( \ell_2^2\mathcal{NN}_n(x_k; \gamma)\right)}{\partial \gamma_{ij}} \right\|
\end{equation}
where $\ell_2^2$ denotes the squared $\ell_2$ norm of the output of the neural network $\mathcal{NN}_n$ applied at $x_k$. The loss function in eq. \ref{eq:loss_MF} is then modified to read: 
\begin{align}
    \mathcal{L}_{MF, MAS}(\gamma^{n}) = &\lambda_{bc} \mathcal{L}_{bc}(\gamma^{n}) + \lambda_{ic} \mathcal{L}_{ic}(\gamma^{n}) + \lambda_{r} \mathcal{L}_{r}(\gamma^{n}) + \lambda_{data} \mathcal{L}_{data}(\gamma^{n}) + \lambda \sum_{i, j} (\gamma_{nl, ij}^n)^2 \nonumber \\
    &+ \lambda_{MAS} \sum_{i, j} \Omega^{n-1}_{ij}\left(\gamma^n_{ij}-\gamma^{n-1}_{ij}\right)^2 \label{eq:loss_MFMAS}
\end{align}

When applying MAS to multifidelity neural networks, we calculate the MAS terms separately: 
\begin{equation}
    \Omega^{n, nl}_{ij} = \frac{1}{N}\sum_{k=1}^N \left\| \frac{\partial\left( \ell_2^2\mathcal{NN}_n^{nl}(x_k; \gamma)\right)}{\partial \gamma_{ij}^{nl}} \right\|, \; \; \; \; 
        \Omega^{n, l}_{ij} = \frac{1}{N}\sum_{k=1}^N \left\| \frac{\partial\left( \ell_2^2\mathcal{NN}_n^{l}(x_k; \gamma)\right)}{\partial \gamma_{ij}^{l}} \right\|
\end{equation}
where $nl$ denotes the nonlinear network and $l$ denotes the linear network. In this way, roughly, the importance in the weights in calculating the linear and nonlinear terms is found separately, instead of determining the  importance in the overall output of the sum of the networks. The parameter $\lambda_{MAS}$ is kept the same for the linear and nonlinear parts.

\subsection{Replay} 
 
In replay, a selection of points in the previously trained domains, $\cup_{i=1}^{n-1} \Omega_{i}$ are selected at each iteration and the residual loss, $\mathcal{L}_r(\gamma^n)$ is evaluated at the points. In this way, the multididelity training still satisfies the PDE across the earlier trained domains. For PINNs, the replay approach only requires knowledge of the geometry of $\cup_{i=1}^{n-1} \Omega_{i}$, and not the value of the output of the model on this domain. 

\subsection{Transfer learning} 
In all cases in this work, the values of the trainable parameters in each subsequent network $\mathcal{NN}_i$, $i \geq 2$, is initialized from the final values of the trainable parameters in the previous network, 
$\mathcal{NN}_{i-1}$. In notation, $\gamma_i^0 = \gamma_{i-1}$. This approach allows for faster training because the network is not initialized randomly. We note that some previous work has found less forgetting by initializing each subsequent network randomly \cite{benzing2020unifying}, and leave the exploration of this option for future work. 

\section{Physics-informed training} \label{sec:pinns}
In this section, we give examples of applying the multifidelity continual learning for physics-informed neural networks in cases where PINNs fail to train. We show that using continual learning in time can improve the accuracy of training a PINN for long-time integration problems, where a single PINN is not sufficient. All hyperparameters used in training are given in Appendix \ref{sec:App}.

\subsection{Pendulum dynamics}\label{sec:pen}
In this section, we consider the gravity pendulum with damping from \cite{wang2023long}. The system is governed by an ODE for $t \in [0, T]$ 
\begin{align}
    \frac{d s_1}{dt} &= s_2, \label{eq:pendulum_1}\\
    \frac{d s_2}{dt} &= -\frac{b}{m} s_2 - \frac{g}{L} \sin(s_1). \label{eq:pendulum_2}
\end{align}
The initial conditions are give by $s_1(0) = s_2(0) = 1$. We take $m=L=1$, $b=0.05$, and $g=9.81$, and we take $T=10$.

\begin{figure}[h]
    \centering
    \includegraphics[width=0.8\textwidth]{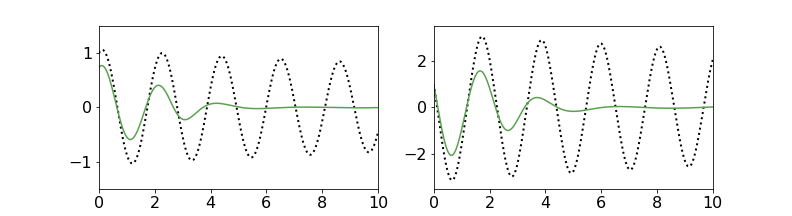}
    \caption{Results from training a single PINN to satisfy Eqs. \ref{eq:pendulum_1} and \ref{eq:pendulum_2} (solid lines) compared with the exact solution (dotted line) for $s_1$ (left) and $s_2$ (right). The results decay to zero quickly and the learned solution does not agree well with the exact solution.}
    \label{fig:PINN_pendulum}
\end{figure}

We first consider a single PINN trained in $t \in [0, 10]$ in Fig. \ref{fig:PINN_pendulum}. The solution quickly goes to zero, showing that a single PINN cannot capture the longtime dynamics of even this simple system. Similar results were shown in \cite{wang2023long}. We will note that there are recent advances that have been developed for improving the training of PINNs for long-time integration problems \cite{wang2023long, wang2022and, meng2020ppinn}. %TO DO: add more citations here. 
In this section, we will explore how continual learning can also allow for accurate solutions over long times by dividing the time domains into subdomains. 

We divide the domain into five subdomains, $\Omega_i = [2(i-1), 2i]$ and train on each domain using both traditional single fidelity continual learning and MF-CL, and SF-CL and MF-CL approaches augmented by replay and MAS. For each case, we calculate the root mean square error (RMSE) of the final output $\mathcal{NN}_5$ on the full domain, $\Omega = [0, 10]$ by 
\begin{equation}
    RMSE = \sqrt{\frac{1}{N}\sum_{j=1}^N \left[\mathcal{NN}_5(t_j)-\mathbf{s}(t_j)\right]^2},
\end{equation}
where $\mathbf{s}$ denotes the exact solution. 
If forgetting is limited, the final solution should have a small RMSE on the full domain.

\begin{figure}[ht!]
    \centering
     \begin{subfigure}[b]{\textwidth}
     \centering
        \includegraphics[width=0.8\textwidth]{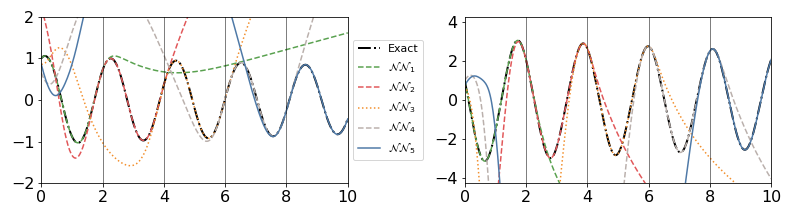}
    \caption{Single fidelity}
    \end{subfigure}
         \begin{subfigure}[b]{\textwidth}
     \centering
        \includegraphics[width=0.8\textwidth]{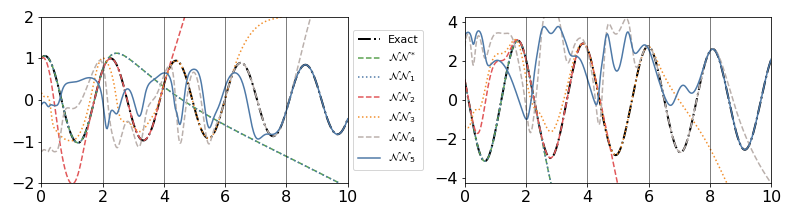}
    \caption{Multifidelity}
    \end{subfigure}
    \caption{Results from training the single fidelity (a) and multifidelity (b) alone  to satisfy Eqs. \ref{eq:pendulum_1} and \ref{eq:pendulum_2} compared with the exact solution (dash-dotted line) for $s_1$ (left) and $s_2$ (right). Of particular importance is the final network, $\mathcal{NN}_5$ (blue solid line), which is trained on $\Omega_5 = [8, 10]$. While the multifidelity results in (b) have significant errors, the are substantially better than the single fidelity results in (a). In the single fidelity training, each network $\mathcal{NN}_i$ is only accurate on the subdomain $\Omega_i$, and extrapolation outside $\Omega_i$ presents significant difficulties. }
    \label{fig:PINN_pendulum_MF_alone}
\end{figure}

\begin{figure}[ht!]
    \centering
     \begin{subfigure}[b]{\textwidth}
     \centering
        \includegraphics[width=0.8\textwidth]{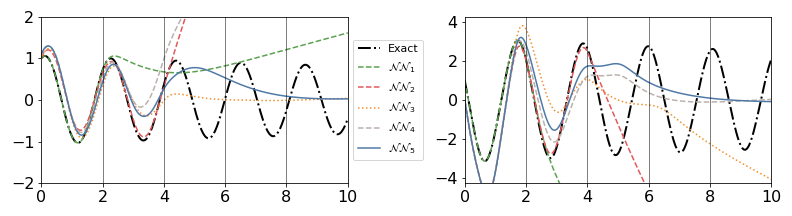}
    \caption{Single fidelity}
    \end{subfigure}
         \begin{subfigure}[b]{\textwidth}
     \centering
        \includegraphics[width=0.8\textwidth]{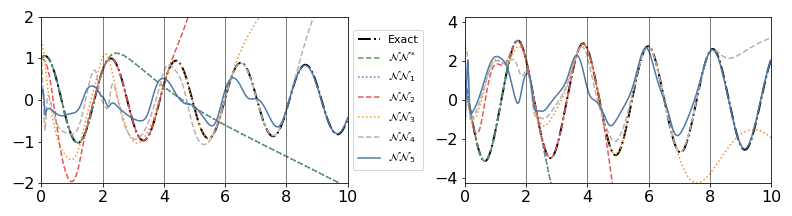}
    \caption{Multifidelity}
    \end{subfigure}
    \caption{Results from training the single fidelity (a) and multifidelity (b) with MAS  to satisfy Eqs. \ref{eq:pendulum_1} and \ref{eq:pendulum_2} compared with the exact solution (dash-dotted line) for $s_1$ (left) and $s_2$ (right). Of particular importance is the final network, $\mathcal{NN}_5$ (blue solid line), which is trained on $\Omega_5 = [8, 10]$. These simulations plotted here have the smallest RMSEs of $\mathcal{NN}_5$ on $\Omega$ of any of the sets of hyperparameters tested. In the single fidelity case, MAS appears to cause restrictions in training that are too strict, and later networks $\mathcal{NN}_i$ are no longer accurate on their respective domains $\Omega_i$. For the multifidelity training, the solutions are accurate across a wider portion of the full domain, and the RMSE is decreased compared with multifidelity training alone.  }
    \label{fig:PINN_pendulum_MF_MAS}
\end{figure}

\begin{table}
\begin{center}
\begin{tabular}{ c|c|c } 
 \hline
  & Single fidelity & Multifidelity \\ 
  \hline
Applied alone & 16.774 & 2.147 \\ 
 Replay & 0.041 & 0.079 \\ 
  MAS  & 1.390 & 1.146 \\
 \hline
\end{tabular}
\caption{RMSE of the final output $\mathcal{NN}_5$ on the full domain for the pendulum problem. For the MAS cases, the network is trained for six values of $\lambda_{MAS}$, and the case with the lowest RMSE is shown in the table above. The replay results have $N=100$ neurons in each hidden layer, see Table \ref{table:3_2} for cases with varying neurons in each hidden layer. }
\label{table:3_1}
\end{center}
\end{table}

It is clear from Table \ref{table:3_1} that replay performs the best in both cases, and significantly better than any other approach. It is no surprise that the SF applied alone case has a large RMSE, as it does not have any incorporation of techniques to limit forgetting. This case is shown in Fig. \ref{fig:PINN_pendulum_MF_alone}a. 

Fig. \ref{fig:PINN_pendulum_MF_MAS} gives the best MAS results for each of the sets of hyperparameters considered, with $\lambda_{MAS} = 100$ for single fidelity and $\lambda_{MAS} = 0.001$ for multifidelity. As is unsurprising given the smaller RMSE, the multifidelity outperforms the single fidelity training with MAS.

\begin{figure}[ht!]
    \centering
     \begin{subfigure}[b]{\textwidth}
     \centering
        \includegraphics[width=0.8\textwidth]{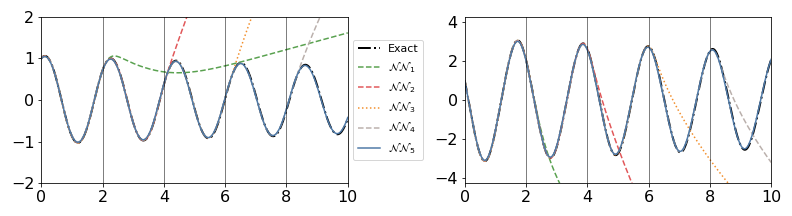}
    \caption{Single fidelity}
    \end{subfigure}
         \begin{subfigure}[b]{\textwidth}
     \centering
        \includegraphics[width=0.8\textwidth]{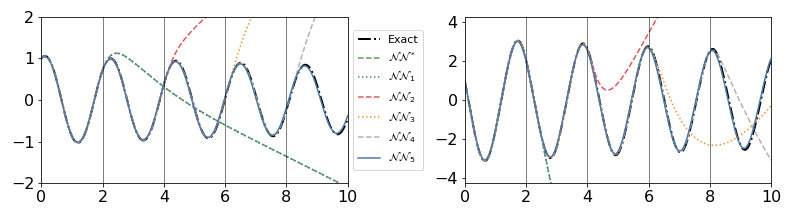}
    \caption{Multifidelity}
    \end{subfigure}
    \caption{Results from training the single fidelity (a) and multifidelity (b) with Replay  to satisfy Eqs. \ref{eq:pendulum_1} and \ref{eq:pendulum_2} compared with the exact solution (dash-dotted line) for $s_1$ (left) and $s_2$ (right). Both cases show very limited forgetting. }
    \label{fig:PINN_pendulum_MF_replay}
\end{figure}

As shown in Fig. \ref{fig:PINN_pendulum_MF_replay}, the SF-replay case does appear to outperform the MF-replay case. However, it is interesting to look at the RMSE as we change the network size in Table \ref{table:3_2}. While the MF-replay case is robust to changes in the network size, the single fidelity case only achieves a small RMSE with a very specific architecture.

\begin{table}
\begin{center}
\begin{tabular}{ c|c|c } 
 \hline
$N$  & Single fidelity & Multifidelity \\ 
  \hline
25 & 3.130 & 0.044 \\ 
50 & 0.107 & 0.027 \\ 
100 & 0.041 & 0.079 \\ 
150 & 0.055 & 0.054 \\ 
200 & 0.091 & 0.038\\ 
 \hline
\end{tabular}
\caption{RMSE of the final output $\mathcal{NN}_5$ on the full domain for the pendulum problem. The SF case has five hidden layers with $N$ neurons each. In the MF case, each nonlinear network has five hidden layers with $N$ neurons. The multifidelity linear network has one hidden layer with 20 neurons.}
\label{table:3_2}
\end{center}
\end{table}

\subsection{Allen-Cahn equation} \label{sec:AC}

\begin{figure}[ht!]
    \centering
    \includegraphics[width=\columnwidth]{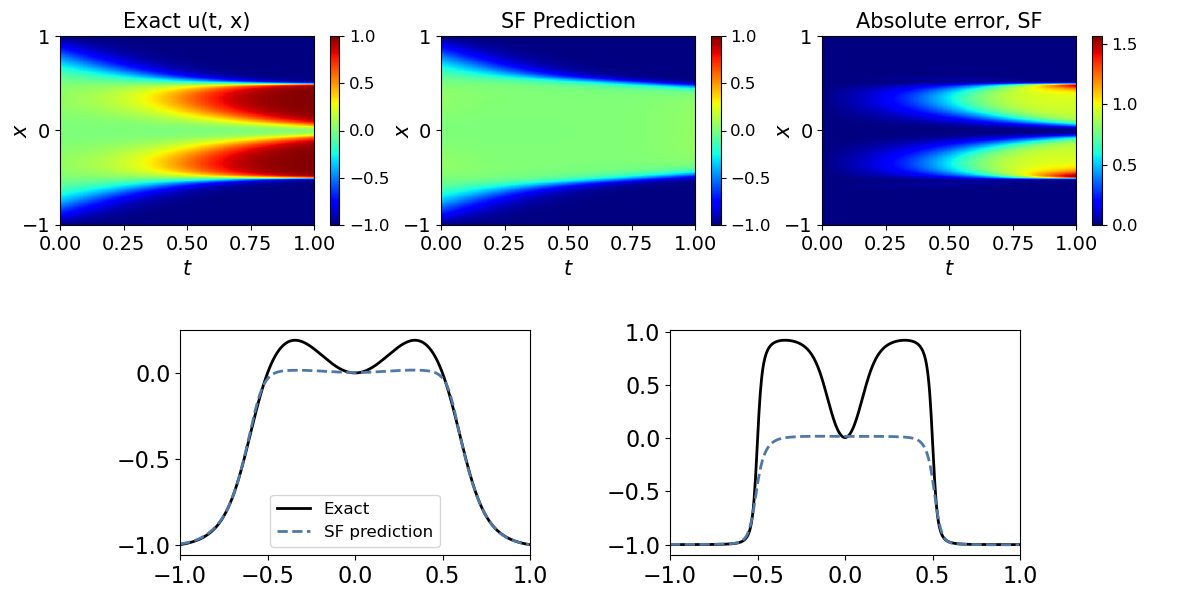}
    \caption{Results from training a single PINN training for the Allen-Cahn equation. The bottom figures are taken at $t = 0.25$ (left) and $t = 0.75$ (right). While the PINN trains well until about 0.3, the solution degrades with increasing $t$. }
    \label{fig:AC_PINN}
\end{figure}

The Allen-Cahn equation is given by 
\begin{align}
    &u_t -c_1^2u_{xx} +5u^3-5u = 0, \;\;\; t\in(0, 1], x\in[-1, 1] \\
    &u(x, 0) = x^2 \cos(\pi x), \;\;\; x\in[-1, 1] \\
    &u(x, t) = u(-x, t), \;\;\; t\in[0, 1], x=-1, x=1 \\
    &u_x(x, t) = u_x(-x, t),\;\;\; t\in[0, 1], x=-1, x=1 
\end{align}
We take $c_1^2 = 0.0001$. The Allen-Cahn equation is notoriously difficult for PINNs to solve by direct application \cite{wight2020solving, rohrhofer2022role}, see Fig. \ref{fig:AC_PINN}. Modifications of PINNs have successfully been able to solve the Allen-Cahn equation, including by using a discrete Runge-Kutta neural network \cite{raissi2019physics}, adaptive sampling of the collocation points \cite{wight2020solving}, and backward compatible PINNs \cite{mattey2022novel}. In this section we show that we can accurately learn the solution to the Allen-Cahn equation by applying the multifidelity continual learning framework.

\begin{figure}[ht!]
    \centering
    \includegraphics[width=\columnwidth]{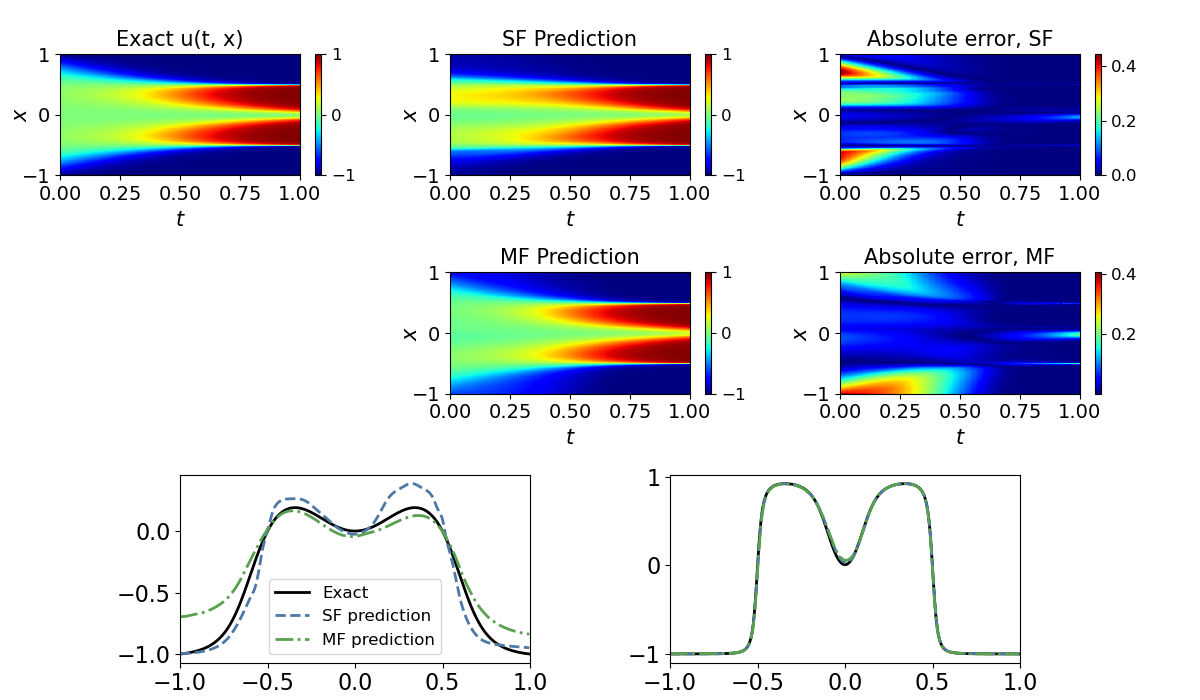}
    \caption{$\mathcal{NN}_4$ results from training a single fidelity and multifidelity PINN training alone for the Allen-Cahn equation. The bottom figures are taken at $t = 0.25$ (left) and $t = 0.75$ (right). The multifidelity results have errors about half as large as those of the single fidelity results.}
    \label{fig:AC_PINN_alone}
\end{figure}

\begin{figure}[ht!]
    \centering
    \includegraphics[width=\columnwidth]{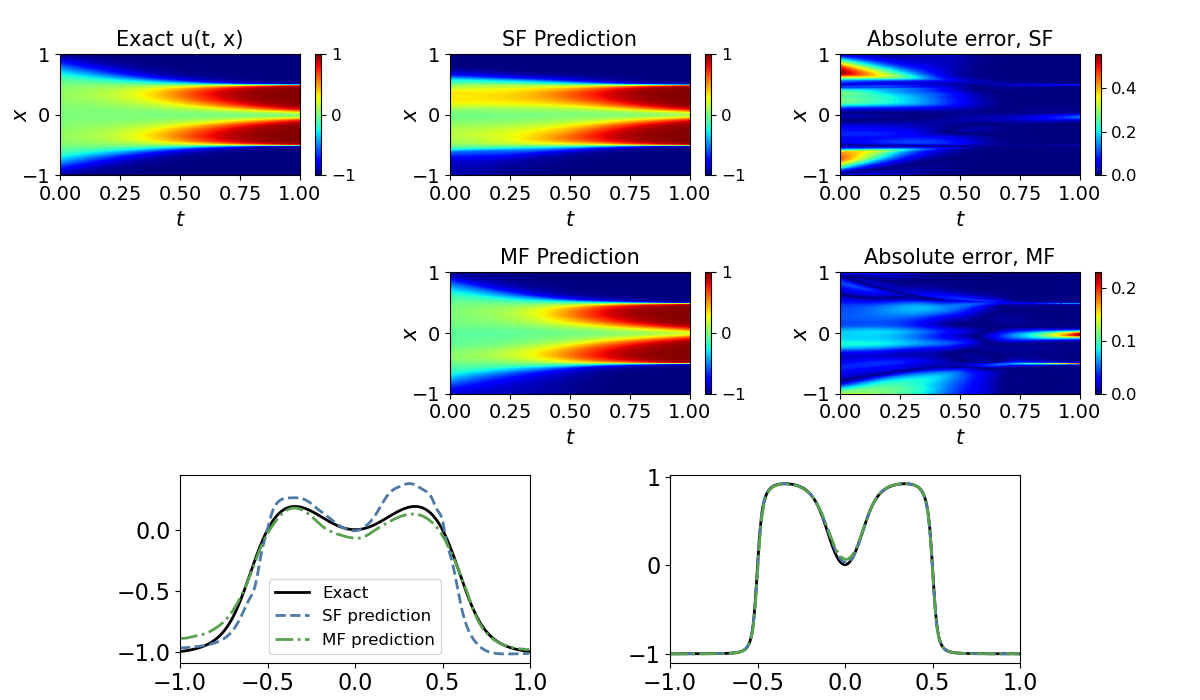}
    \caption{$\mathcal{NN}_4$ results from training a single fidelity and multifidelity PINN training with MAS for the Allen-Cahn equation. The bottom figures are taken at $t = 0.25$ (left) and $t = 0.75$ (right). These results represent the best MAS results from all sets of hyperparameters considered. The multifidelity results have errors about a quarter as large as those of the single fidelity results.}
    \label{fig:AC_PINN_MAS}
\end{figure}

\begin{figure}[ht!]
    \centering
    \includegraphics[width=\columnwidth]{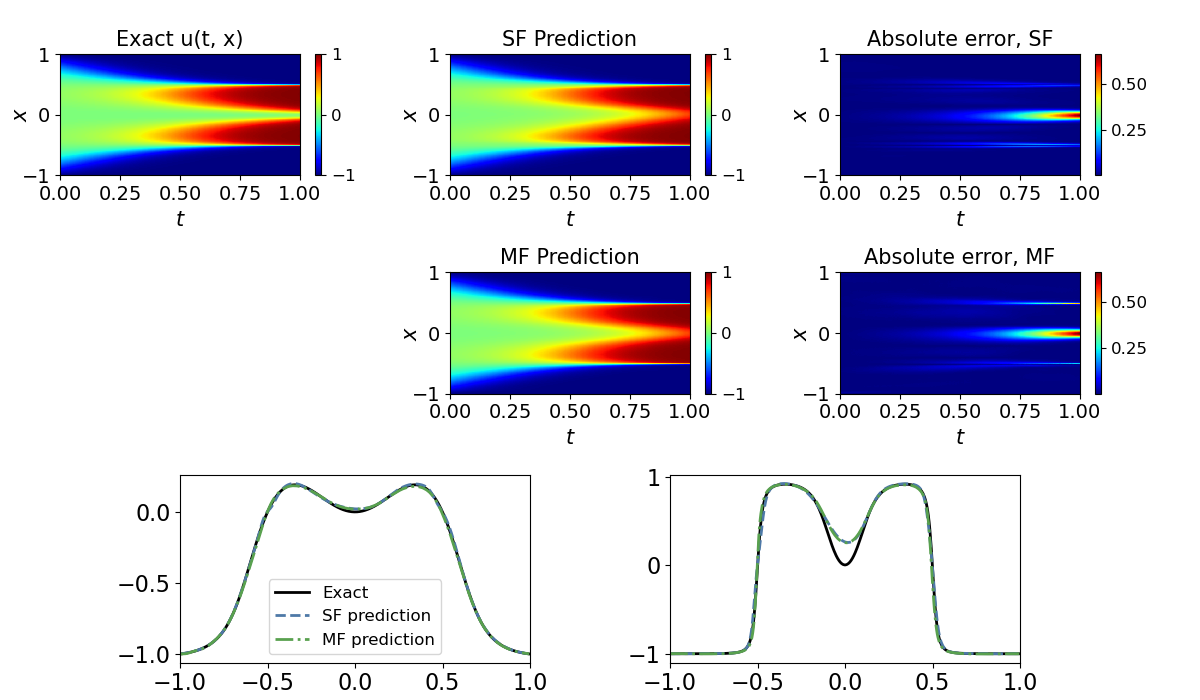}
    \caption{$\mathcal{NN}_4$ results from training a single fidelity and multifidelity PINN training with replay for the Allen-Cahn equation. The bottom figures are taken at $t = 0.25$ (left) and $t = 0.75$ (right).}
    \label{fig:AC_PINN_replay}
\end{figure}

\begin{table}
\begin{center}
\begin{tabular}{ c|c|c } 
 \hline
  & Single fidelity & Multifidelity \\ 
  \hline
 No CL & 0.126 & 0.128 \\ 
 Replay & 0.087 & 0.084 \\ 
  MAS  & 0.146 & 0.056 \\
 \hline
\end{tabular}
\caption{Relative RMSE of the final output $\mathcal{NN}_4$ on the full domain for the Allen-Cahn equation. For the MAS cases, the network is trained for seven values of $\lambda_{MAS}$, and the case with the lowest RMSE is shown in the table above.}
\label{table:AC_error}
\end{center}
\end{table}

We divide the domain into four subdomains, $\Omega_i = [2(i-1), 2i]$, and report the relative RMSE of $\mathcal{NN}_4$ on the full domain $\Omega.$ When the multifidelity and single fidelity methods are trained alone, in Fig. \ref{fig:AC_PINN_alone}, they have approximately equal relative RMSEs. MAS and replay both improve the results, in Figs. \ref{fig:AC_PINN_MAS} and \ref{fig:AC_PINN_replay}, respectively. A summary of the results is given in Table \ref{table:AC_error}.

\section{Data-informed training} \label{sec:data}

\subsection{Batteries} \label{sec:batteries}
This is a case where if an additional dataset is added, it not is clear a priori which subdomain it lies in. Therefore, it is essential that the final model can predict the current accurately for the entire domain without forgetting. 

\begin{figure}[ht!]
    \centering
    \includegraphics[width=1\columnwidth]{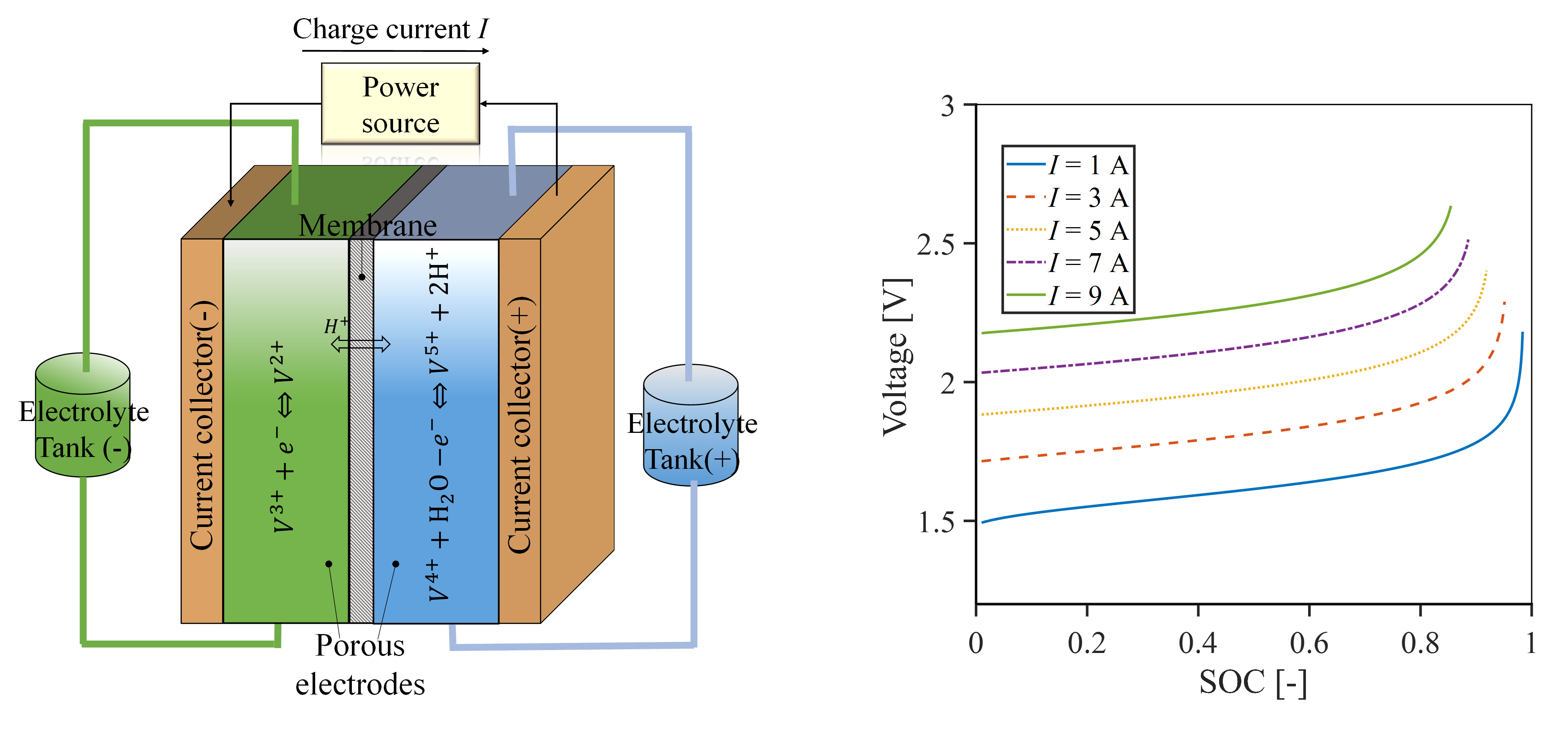}
    \caption{The VRFB system used for battery data generation (left).  Sample charge curve distribution at different charge current (right).}
    \label{fig:Battery_Model}
\end{figure}

For testing, a vanadium redox-flow battery (VRFB) system was selected to generate datasets. The left image in Fig. \ref{fig:Battery_Model} shows a typical configuration of a VRFB, which consists of electrodes, current collectors and a membrane separator. The negative and positive side have a storage tank each to store the redox couple of $\text{V}^{2+}/\text{V}^{3+}$ and $\text{V}^{4+}/\text{V}^{5+}$, respectively. We applied the MFCL method for the problem of identifying the applied charge current from a given charge voltage curve. To generate the VRFB charge curve dataset, a highly computationally efficient 2-D analytical model was utilized \cite{CHEN_2D, CHEN_Analytical}. This model fully resolves the coupled physics of active species transport, electrochemical reaction kinetics, and fluid dynamics within the battery cell, thereby providing a faithful representation of the VRFB system. Further details on the model and its parameters can be found in \cite{CHEN_2D}. Typical charge curves are visualized in the right plot of Fig. \ref{fig:Battery_Model} for five selected current levels. For a given charge current, the battery voltage ($E$) is calculated at different state-of-charge (SOC) values to form the charge curve which is used as input data. The applied charge current $I$ which gives rise to the charge curve is the output quantity we want to predict.

\begin{figure}[h]
    \centering
         \begin{subfigure}[b]{\textwidth}
\centering
    \includegraphics[width=0.7\columnwidth]{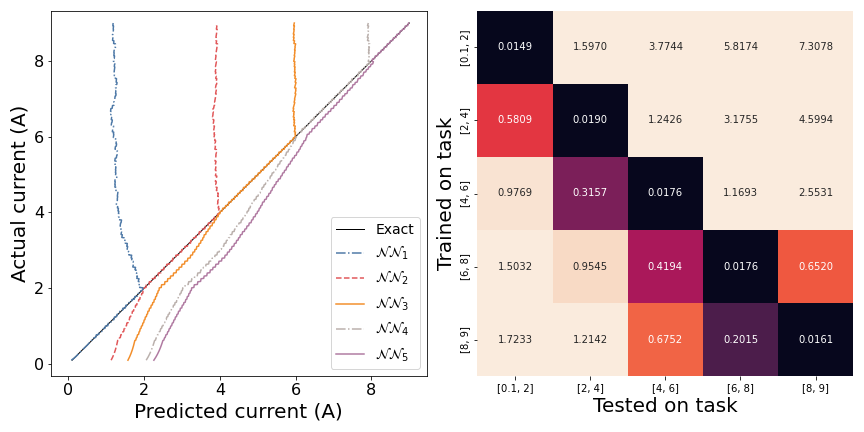}
    \caption{Single fidelity}
    \end{subfigure}
             \begin{subfigure}[b]{\textwidth}
\centering
    \includegraphics[width=0.7\columnwidth]{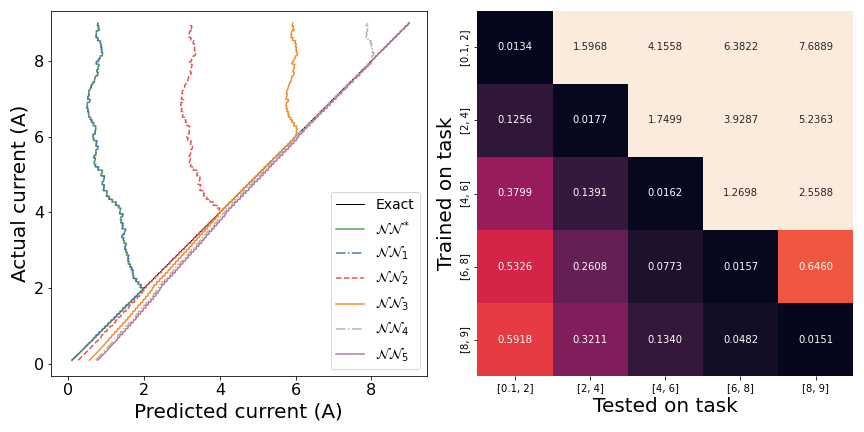}
        \caption{Multifidelity}
    \end{subfigure}
    \caption{Results from the single fidelity (a) and multifidelity (b) training alone for the battery test case. The left column has the network outputs of each task on all the tasks, and the right column shows the RMSE of each task tested on each other task. The results in this figure use the narrow architecture, with three hidden layers and 40 neurons per layer. }
    \label{fig:Battery_alone}
\end{figure}

We divide the data set in five sets by charge current and train with and without MAS in the single fidelity case. The subdomains are $\Omega_1=[0.1, 2)$, $\Omega_2=[2, 4)$, $\Omega_3=[4, 6)$, $\Omega_4=[6, 8)$, and $\Omega_5=[8, 9]$. The errors are calculated by the RMSE of the output of $\mathcal{NN}_5$ on a test test selected from $\Omega = \cup_{i=1}^5 \Omega_i$. We test two network architectures, a wide network which has two hidden layers with 80 neurons each, and a deeper and narrower network which has three hidden layers with 40 neurons each. We first train with the single fidelity and multifdelity approaches alone, see Fig. \ref{fig:Battery_alone}. The multifidelity continual learning results show less forgetting than those from the single fidelity continual learning. 

\begin{figure}[h]
    \centering
         \begin{subfigure}[b]{\textwidth}
\centering
    \includegraphics[width=0.7\columnwidth]{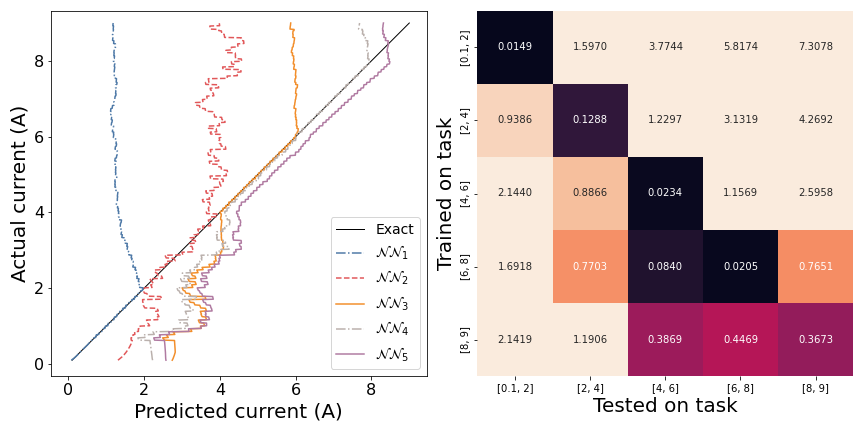}
    \caption{Single fidelity}
    \end{subfigure}
             \begin{subfigure}[b]{\textwidth}
\centering
    \includegraphics[width=0.7\columnwidth]{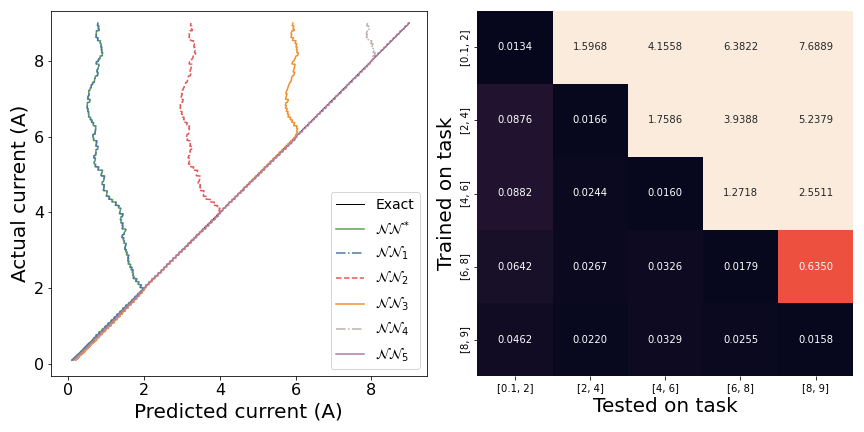}
        \caption{Multifidelity}
    \end{subfigure}
    \caption{Results from the single fidelity (a) and multifidelity (b) training with MAS for the battery test case. The single fidelity case struggles to train accurately, while multifidelity has very limited forgetting. The left column has the network outputs of each task on all the tasks, and the right column shows the RMSE of each task tested on each other task. The results shown represent the best output from the MAS hyperparameters tested. For the single fidelity case, the results are from the narrow network with $\lambda_{MAS} = 100$, and for the multifidelity case, the results are from the wide network with $\lambda_{MAS} = 0.001$. }
    \label{fig:Battery_MAS}
\end{figure}

\begin{figure}[h]
    \centering

    \includegraphics[width=0.4\columnwidth]{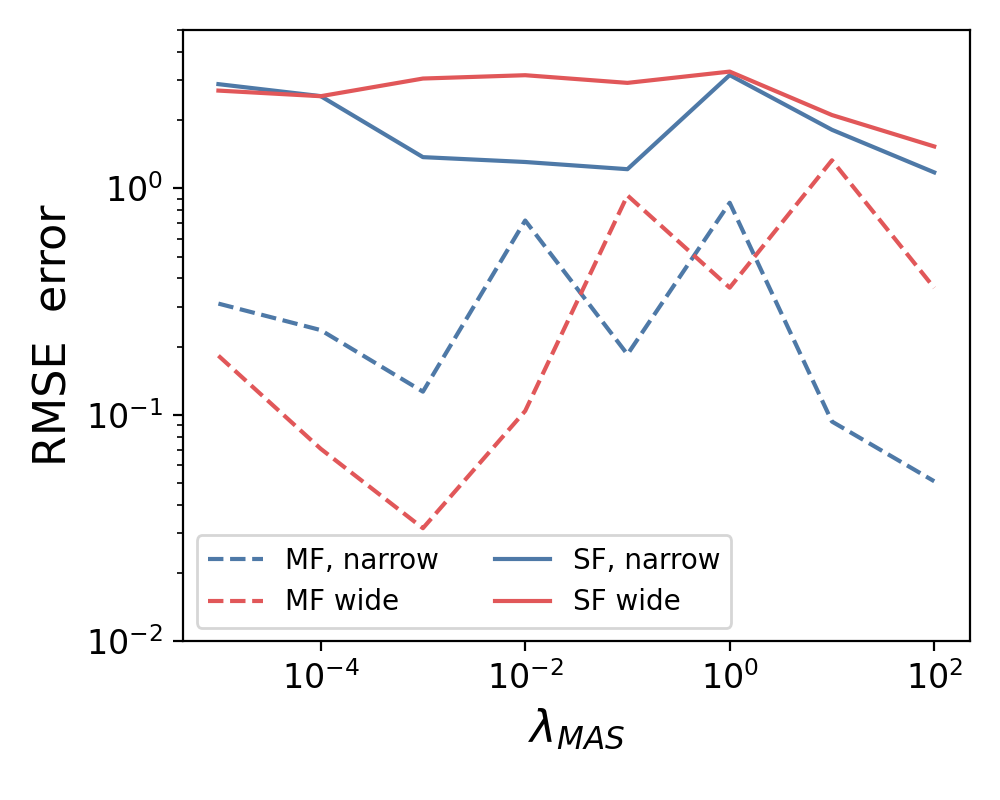}
    \caption{Comparison of the RMSE with MAS for the redox flow battery test case. The RMSE is lower for almost all the multifidelity test cases in comparison with the single fidelity test cases. }
    \label{fig:Battery_MAS_comp}
\end{figure}

We then consider the impact of adding MAS. We consider the narrow and wide networks with and without MAS scaling, for a total of four cases. The multifidelity MAS results show significant improvement, see Fig. \ref{fig:Battery_MAS}. In Fig. \ref{fig:Battery_MAS_comp}, we compare the performance across the value of the MAS hyperparameter $\lambda_{MAS}$. We see that the single fidelity approach performance is robust, since it is insensitive to the value of $\lambda_{MAS}.$ However, it is not very accurate. On the other hand, the multifidelity approach can be substantially more accurate than the single fidelity approach for most values of $\lambda_{MAS}.$ Overall, the multifidelity results significantly outperform the single fidelity results.

\subsection{Energy consumption} \label{sec:energy}

To provide a second example of data-informed continual learning, we consider the city-scale daily energy consumption dataset from \cite{wang2021predicting}. The dataset consists of daily energy usage for three metropolitan areas, New York, Sacramento, and Los Angeles, along with daily weather data. Three years of data are used as a test set, with an additional year as a test set. 

\begin{figure}[h]
    \centering

    \includegraphics[width=0.8\columnwidth]{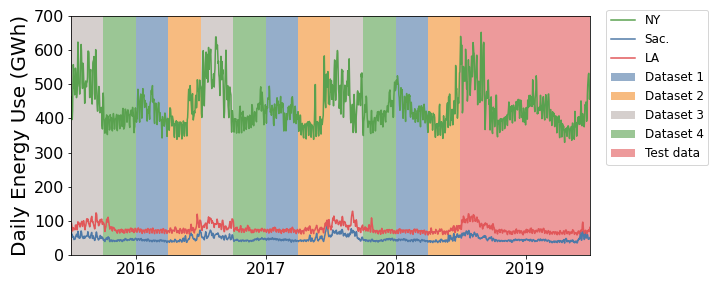}
    \caption{Illustration of the datasets used in the energy consumption example. }
    \label{fig:Energy_data}
\end{figure}

Energy usage depends strongly on the weather, with air conditioner usage in the warmer months and heating in the winter months. Therefore, to provide different tasks to the continual learning training, we divide the three years of training data by quarter. Task 1 has training data from January to March, Task 2 has training data from April to June, Task 3 has training data from July to September, and Task 4 has training data from October to December. The test set for all tasks is to predict the energy usage from July 2018 to June 2019. An illustration of the testing and training data divided into tasks is given in Fig. \ref{fig:Energy_data}.

\begin{figure}[h]
    \centering

    \includegraphics[width=0.8\columnwidth]{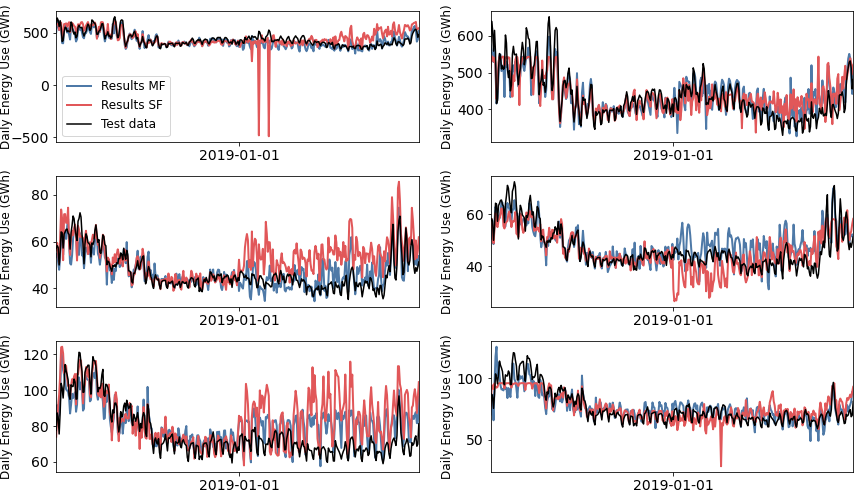}
    \caption{Results of the energy consumption problem. (Left) results without MAS. (Right) results with MAS.}
    \label{fig:Energy_results}
\end{figure}

\begin{figure}[h]
    \centering

    \includegraphics[width=0.6\columnwidth]{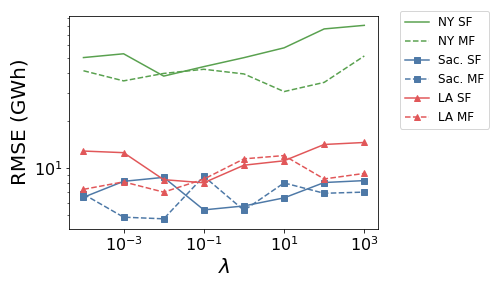}
    \caption{Comparison of the RMSEs generated by training with MAS for the energy consumption problem.}
    \label{fig:Energy_MAS_comp}
\end{figure}

We train both single fidelity and multifidelity networks with and without MAS. We consider a range of $\lambda \in [0.001, 100]$. We also compare with training a network without continual learning. In this case, a single fidelity DNN receives \emph{all} of the training data from all four tasks, to try and predict the energy usage from July 2018 to June 2019. This case serves as a benchmark for the reasonable level of error we can expect from our model using continual learning. 
% as shown in Fig. \ref{fig:Energy_results_noCL}.
The results are shown in Table \ref{table:Energy_error}. We note that in all cases, the multifidelity continual learning approach outperforms the single fidelity continual learning. Including MAS does improve the results, as shown in Fig. \ref{fig:Energy_results}. The continual learning methods do perform worse than the case with no continual learning, which is expected because they never have access to all the training data simultaneously. A comparison of the RMSE for each value of $\lambda_{MAS}$ tested is given in Fig. \ref{fig:Energy_MAS_comp}. We note that overall, the multifidelity approach with MAS is more robust than the single fidelity training with MAS, resulting in a smaller RMSE across a range of $\lambda_{MAS}.$ 

\begin{table}
\begin{center}
\begin{tabular}{ c||  c | c|c | c| c} 
 \hline
 City & No CL & Single fidelity & Multifidelity & Single fidelity, MAS  & Multifidelity, MAS\\ 
  \hline
New York &  16.6289 &    95.8546 &    36.8522 &    38.2728 &    30.5356 \\
Sacramento &   2.6952 &     8.9388 &     4.4678 &     5.4253 &     4.7614 \\
Los Angeles &    5.5859 &    15.1774 &     9.8254 &     8.0617 &     7.0310 \\
 \hline
\end{tabular}
\caption{RMSE (GWh) of the final output $\mathcal{NN}_4$ on the full test domain for the energy consumption case. For the MAS cases, the network is trained for seven values of $\lambda_{MAS}$, and the case with the lowest RMSE is shown in the table above.}
\label{table:Energy_error}
\end{center}
\end{table}

.

%\begin{figure}[h]
 %   \centering

   % \includegraphics[width=0.48\columnwidth]{Final_figs/Figure18.png}
  %  \caption{Results of the energy consumption problem without continual learning}
    %\label{fig:Energy_results_noCL}
%\end{figure}

\section{Discussion and future work}

We have introduced a novel continual learning method based on multifidelity deep neural networks. The premise of the method is the existence of correlations between the output of previously trained models and the desired output of the model on the current training dataset. The discovery and use of these correlations can  limit catastrophic forgetting. On its own, the multifidelity continual learning method has shown robustness and limited forgetting across several datasets for physics-informed and data-driven training examples. Additionally, it can be combined with existing continual learning methods, including replay and memory aware synapses (MAS), to further limit catastrophic forgetting. 

The proposed continual learning method is especially suited for physical problems where the data satisfy the same physical laws on each domain, or for a physics-informed neural network, because in these cases we expect there to be a strong correlation between the output of the previous model and the model on the current training domain. As a result of exploiting the correlation between data in the various domains instead of training from scratch for each domain, the method can afford to continue learning in new domains using smaller networks. Specifically, its training accuracy is more robust to the size of the network employed in the new domain. This can lead to computational savings during both training and inference. The approach is particularly suited for situations where privacy concerns can limit access to prior datasets. It can also offer new possibilities in the area of federated learning by allowing the design of new algorithms for processing sensor data in a distributed fashion. These topics are under investigation and results will be reported in a future publication.  

\section{Data availability }
All code and data needed to reproduce these results will be made available at \url{https://github.com/pnnl/Multifidelity_continual_learning/}.

\section{Acknowledgements}
This research was supported by the Energy Storage Materials Initiative (ESMI), under the Laboratory Directed Research and Development (LDRD) Program at Pacific Northwest National Laboratory (PNNL).  The computational work was performed using PNNL Institutional Computing at Pacific Northwest National Laboratory. PNNL is a multi-program national laboratory operated for the U.S. Department of Energy (DOE) by Battelle Memorial Institute under Contract No. DE-AC05-76RL01830.

\bibliographystyle{unsrt}  
\bibliography{references}  

\begin{thebibliography}{10}

\bibitem{parisi2019continual}
German~I Parisi, Ronald Kemker, Jose~L Part, Christopher Kanan, and Stefan
  Wermter.
\newblock Continual lifelong learning with neural networks: A review.
\newblock {\em Neural networks}, 113:54--71, 2019.

\bibitem{verwimp2021rehearsal}
Eli Verwimp, Matthias De~Lange, and Tinne Tuytelaars.
\newblock Rehearsal revealed: The limits and merits of revisiting samples in
  continual learning.
\newblock In {\em Proceedings of the IEEE/CVF International Conference on
  Computer Vision}, pages 9385--9394, 2021.

\bibitem{zenke2017continual}
Friedemann Zenke, Ben Poole, and Surya Ganguli.
\newblock Continual learning through synaptic intelligence.
\newblock In {\em International conference on machine learning}, pages
  3987--3995. PMLR, 2017.

\bibitem{kirkpatrick2017overcoming}
James Kirkpatrick, Razvan Pascanu, Neil Rabinowitz, Joel Veness, Guillaume
  Desjardins, Andrei~A Rusu, Kieran Milan, John Quan, Tiago Ramalho, Agnieszka
  Grabska-Barwinska, et~al.
\newblock Overcoming catastrophic forgetting in neural networks.
\newblock {\em Proceedings of the National Academy of Sciences},
  114(13):3521--3526, 2017.

\bibitem{aljundi2018memory}
Rahaf Aljundi, Francesca Babiloni, Mohamed Elhoseiny, Marcus Rohrbach, and
  Tinne Tuytelaars.
\newblock Memory aware synapses: Learning what (not) to forget.
\newblock In {\em Proceedings of the European Conference on Computer Vision
  (ECCV)}, pages 139--154, 2018.

\bibitem{de2021continual}
Matthias De~Lange, Rahaf Aljundi, Marc Masana, Sarah Parisot, Xu~Jia,
  Ale{\v{s}} Leonardis, Gregory Slabaugh, and Tinne Tuytelaars.
\newblock A continual learning survey: Defying forgetting in classification
  tasks.
\newblock {\em IEEE transactions on pattern analysis and machine intelligence},
  44(7):3366--3385, 2021.

\bibitem{hsu2018re}
Yen-Chang Hsu, Yen-Cheng Liu, Anita Ramasamy, and Zsolt Kira.
\newblock Re-evaluating continual learning scenarios: A categorization and case
  for strong baselines.
\newblock {\em arXiv preprint arXiv:1810.12488}, 2018.

\bibitem{rusu2016progressive}
Andrei~A Rusu, Neil~C Rabinowitz, Guillaume Desjardins, Hubert Soyer, James
  Kirkpatrick, Koray Kavukcuoglu, Razvan Pascanu, and Raia Hadsell.
\newblock Progressive neural networks.
\newblock {\em arXiv preprint arXiv:1606.04671}, 2016.

\bibitem{wen2020batchensemble}
Yeming Wen, Dustin Tran, and Jimmy Ba.
\newblock Batchensemble: an alternative approach to efficient ensemble and
  lifelong learning.
\newblock {\em arXiv preprint arXiv:2002.06715}, 2020.

\bibitem{pfeiffer2020adapterhub}
Jonas Pfeiffer, Andreas R{\"u}ckl{\'e}, Clifton Poth, Aishwarya Kamath, Ivan
  Vuli{\'c}, Sebastian Ruder, Kyunghyun Cho, and Iryna Gurevych.
\newblock Adapterhub: A framework for adapting transformers.
\newblock {\em arXiv preprint arXiv:2007.07779}, 2020.

\bibitem{bereska2022continual}
Leonard Bereska and Efstratios Gavves.
\newblock Continual learning of dynamical systems with competitive federated
  reservoir computing.
\newblock In {\em Conference on Lifelong Learning Agents}, pages 335--350.
  PMLR, 2022.

\bibitem{munkhdalai2017meta}
Tsendsuren Munkhdalai and Hong Yu.
\newblock Meta networks.
\newblock In {\em International conference on machine learning}, pages
  2554--2563. PMLR, 2017.

\bibitem{vladymyrov2023continual}
Max Vladymyrov, Andrey Zhmoginov, and Mark Sandler.
\newblock Continual few-shot learning using hypertransformers.
\newblock {\em arXiv preprint arXiv:2301.04584}, 2023.

\bibitem{karniadakis2021physics}
George~Em Karniadakis, Ioannis~G Kevrekidis, Lu~Lu, Paris Perdikaris, Sifan
  Wang, and Liu Yang.
\newblock Physics-informed machine learning.
\newblock {\em Nature Reviews Physics}, 3(6):422--440, 2021.

\bibitem{baker2019workshop}
Nathan Baker, Frank Alexander, Timo Bremer, Aric Hagberg, Yannis Kevrekidis,
  Habib Najm, Manish Parashar, Abani Patra, James Sethian, Stefan Wild, et~al.
\newblock Workshop report on basic research needs for scientific machine
  learning: Core technologies for artificial intelligence.
\newblock Technical report, USDOE Office of Science (SC), Washington, DC
  (United States), 2019.

\bibitem{cuomo2022scientific}
Salvatore Cuomo, Vincenzo~Schiano Di~Cola, Fabio Giampaolo, Gianluigi Rozza,
  Maziar Raissi, and Francesco Piccialli.
\newblock Scientific machine learning through physics--informed neural
  networks: Where we are and what’s next.
\newblock {\em Journal of Scientific Computing}, 92(3):88, 2022.

\bibitem{jin2021nsfnets}
Xiaowei Jin, Shengze Cai, Hui Li, and George~Em Karniadakis.
\newblock Nsfnets (navier-stokes flow nets): Physics-informed neural networks
  for the incompressible navier-stokes equations.
\newblock {\em Journal of Computational Physics}, 426:109951, 2021.

\bibitem{raissi2020hidden}
Maziar Raissi, Alireza Yazdani, and George~Em Karniadakis.
\newblock Hidden fluid mechanics: Learning velocity and pressure fields from
  flow visualizations.
\newblock {\em Science}, 367(6481):1026--1030, 2020.

\bibitem{cai2021physics}
Shengze Cai, Zhiping Mao, Zhicheng Wang, Minglang Yin, and George~Em
  Karniadakis.
\newblock Physics-informed neural networks (pinns) for fluid mechanics: A
  review.
\newblock {\em Acta Mechanica Sinica}, 37(12):1727--1738, 2021.

\bibitem{Joglekar_2023}
Archis~S Joglekar and Alexander G~R Thomas.
\newblock Machine learning of hidden variables in multiscale fluid simulation.
\newblock {\em Machine Learning: Science and Technology}, 4(3):035049, sep
  2023.

\bibitem{liu2019multi}
Dehao Liu and Yan Wang.
\newblock Multi-fidelity physics-constrained neural network and its application
  in materials modeling.
\newblock {\em Journal of Mechanical Design}, 141(12), 2019.

\bibitem{chen2020physics}
Yuyao Chen, Lu~Lu, George~Em Karniadakis, and Luca Dal~Negro.
\newblock Physics-informed neural networks for inverse problems in nano-optics
  and metamaterials.
\newblock {\em Optics express}, 28(8):11618--11633, 2020.

\bibitem{fang2019deep}
Zhiwei Fang and Justin Zhan.
\newblock Deep physical informed neural networks for metamaterial design.
\newblock {\em IEEE Access}, 8:24506--24513, 2019.

\bibitem{mao2020physics}
Zhiping Mao, Ameya~D Jagtap, and George~Em Karniadakis.
\newblock Physics-informed neural networks for high-speed flows.
\newblock {\em Computer Methods in Applied Mechanics and Engineering},
  360:112789, 2020.

\bibitem{misyris2020physics}
George~S Misyris, Andreas Venzke, and Spyros Chatzivasileiadis.
\newblock Physics-informed neural networks for power systems.
\newblock In {\em 2020 IEEE Power \& Energy Society General Meeting (PESGM)},
  pages 1--5. IEEE, 2020.

\bibitem{huang2022applications}
Bin Huang and Jianhui Wang.
\newblock Applications of physics-informed neural networks in power systems-a
  review.
\newblock {\em IEEE Transactions on Power Systems}, 38(1):572--588, 2022.

\bibitem{moya2023dae}
Christian Moya and Guang Lin.
\newblock Dae-pinn: a physics-informed neural network model for simulating
  differential algebraic equations with application to power networks.
\newblock {\em Neural Computing and Applications}, 35(5):3789--3804, 2023.

\bibitem{raissi2019physics}
Maziar Raissi, Paris Perdikaris, and George~E Karniadakis.
\newblock Physics-informed neural networks: A deep learning framework for
  solving forward and inverse problems involving nonlinear partial differential
  equations.
\newblock {\em Journal of Computational Physics}, 378:686--707, 2019.

\bibitem{wang2023long}
Sifan Wang and Paris Perdikaris.
\newblock Long-time integration of parametric evolution equations with
  physics-informed deeponets.
\newblock {\em Journal of Computational Physics}, 475:111855, 2023.

\bibitem{mattey2022novel}
Revanth Mattey and Susanta Ghosh.
\newblock A novel sequential method to train physics informed neural networks
  for allen cahn and cahn hilliard equations.
\newblock {\em Computer Methods in Applied Mechanics and Engineering},
  390:114474, 2022.

\bibitem{dekhovich2023ipinns}
Aleksandr Dekhovich, Marcel~HF Sluiter, David~MJ Tax, and Miguel~A Bessa.
\newblock ipinns: Incremental learning for physics-informed neural networks.
\newblock {\em arXiv preprint arXiv:2304.04854}, 2023.

\bibitem{rasht2022physics}
Majid Rasht-Behesht, Christian Huber, Khemraj Shukla, and George~Em
  Karniadakis.
\newblock Physics-informed neural networks (pinns) for wave propagation and
  full waveform inversions.
\newblock {\em Journal of Geophysical Research: Solid Earth},
  127(5):e2021JB023120, 2022.

\bibitem{meng2020composite}
Xuhui Meng and George~Em Karniadakis.
\newblock A composite neural network that learns from multi-fidelity data:
  Application to function approximation and inverse pde problems.
\newblock {\em Journal of Computational Physics}, 401:109020, 2020.

\bibitem{benzing2020unifying}
Frederik Benzing.
\newblock Unifying regularisation methods for continual learning.
\newblock {\em arXiv preprint arXiv:2006.06357}, 2020.

\bibitem{wang2022and}
Sifan Wang, Xinling Yu, and Paris Perdikaris.
\newblock When and why pinns fail to train: A neural tangent kernel
  perspective.
\newblock {\em Journal of Computational Physics}, 449:110768, 2022.

\bibitem{meng2020ppinn}
Xuhui Meng, Zhen Li, Dongkun Zhang, and George~Em Karniadakis.
\newblock Ppinn: Parareal physics-informed neural network for time-dependent
  pdes.
\newblock {\em Computer Methods in Applied Mechanics and Engineering},
  370:113250, 2020.

\bibitem{wight2020solving}
Colby~L Wight and Jia Zhao.
\newblock Solving {A}llen-{C}ahn and {C}ahn-{H}illiard equations using the
  adaptive physics informed neural networks.
\newblock {\em arXiv preprint arXiv:2007.04542}, 2020.

\bibitem{rohrhofer2022role}
Franz~M Rohrhofer, Stefan Posch, Clemens G{\"o}{\ss}nitzer, and Bernhard~C
  Geiger.
\newblock On the role of fixed points of dynamical systems in training
  physics-informed neural networks.
\newblock {\em arXiv preprint arXiv:2203.13648}, 2022.

\bibitem{CHEN_2D}
Yunxiang Chen, Jie Bao, Zhijie Xu, Peiyuan Gao, Litao Yan, Soowhan Kim, and Wei
  Wang.
\newblock A two-dimensional analytical unit cell model for redox flow battery
  evaluation and optimization.
\newblock {\em Journal of Power Sources}, 506:230192, 2021.

\bibitem{CHEN_Analytical}
Yunxiang Chen, Zhijie Xu, Chao Wang, Jie Bao, Brian Koeppel, Litao Yan, Peiyuan
  Gao, and Wei Wang.
\newblock Analytical modeling for redox flow battery design.
\newblock {\em Journal of Power Sources}, 482:228817, 2021.

\bibitem{wang2021predicting}
Zhe Wang, Tianzhen Hong, Han Li, and Mary~Ann Piette.
\newblock Predicting city-scale daily electricity consumption using data-driven
  models.
\newblock {\em Advances in Applied Energy}, 2:100025, 2021.

\end{thebibliography}

\newpage
\section{Appendix}\label{sec:App}

In this section we report the training parameters used to train the results reported above. 

\begin{landscape}
\begin{table}
    \centering
    \begin{tabular}{l | c | c |c | c}
        & Sec. \ref{sec:pen} & Sec. \ref{sec:AC} &Sec. \ref{sec:batteries}  & Sec. \ref{sec:energy} \\ \hline
       SF $\mathcal{NN}_1$ learning rate  & (1e-3, 2000, 0.95) & (1e-4, 2000, 0.99) & (1e-3, 2000, 0.9) & (1e-3, 2000, 0.99) \\
       SF $\mathcal{NN}_i$, $i > 1$ learning rate & (1e-3, 2000, 0.95) & (1e-4, 2000, 0.99) & (1e-3, 2000, 0.95) & (1e-3, 2000, 0.99) \\
       SF activation function & \texttt{swish} & \texttt{tanh} & \texttt{relu} & \texttt{tanh}\\
       SF architecture & [1, 100, 100, 100, 100, 100, 2]  & [2, 200, 200, 200, 200, 200, 1] & [784, 80, 80, 1] or & [5, 100, 100, 100, 1] \\
               &  &  & [784, 40, 40, 40, 1] & \\
       SF batch size & 100 & 500 & 10 & 100 \\
       SF boundary condition batch size  & 1 & 100 & -- & -- \\
       SF iterations & 50,000 & 100,000 & 100,000 & 100,000 \\
       MF $\mathcal{NN}^*$ learning rate  & (1e-3, 2000, 0.99) & (1e-4, 2000, 0.99) & (1e-3, 2000, 0.9) & (1e-3, 2000, 0.99) \\
       MF $\mathcal{NN}^*$ activation function  & \texttt{swish} & \texttt{tanh}\texttt{tanh} & \texttt{relu} & \texttt{tanh} \\
       MF $\mathcal{NN}^*$ architecture  & [1, 200, 200, 200, 1] & [2, 200, 200, 200, 200, 200, 1] & [784, 40, 40, 40, 1] & [5, 100, 100, 100, 1]\\
       MF $\mathcal{NN}^*$ iterations  & 50,000 & 50,000  & 100,000 & 100,000 \\
       MF batch size  & 200 & 1000 & 10 & 100 \\
       MF boundary condition batch size  & 1 & 100 & -- & -- \\
       MF $\mathcal{NN}_i$ learning rate  & (1e-3, 2000, 0.99)  & (5e-4, 2000, 0.99) & (1e-3, 2000, 0.95) & (1e-3, 2000, 0.99)\\
       MF $\mathcal{NN}_i$ activation function  & \texttt{swish} & \texttt{tanh} & \texttt{relu} & \texttt{tanh} \\
       MF $\mathcal{NN}_i$ nonlinear architecture  & [3, 100, 100, 100, 100, 100, 2] & [3, 200, 200, 200, 200, 200, 1] &  [785, 80, 80, 1] or  & [6, 100, 100, 100, 1] \\
     &  &  &[785, 40, 40, 40, 1]  &  \\
       MF $\mathcal{NN}_i$ linear architecture  & [2, 20, 2] & [1, 20, 1] & [1, 1] & [1, 5, 1]\\
       MF $\mathcal{NN}_1$ iterations  & 50,000 & 50,000 & 100,000 & 100,000 \\
       MF $\mathcal{NN}_i$, $i > 1$ iterations  & 100,000 &  100,000& 100,000 & 100,000 \\
       MF MAS samples  & 1000 & 3000 & 400 & 1201 \\
       SF MAS samples  & 1000 & 3000 & 400 & 1201 \\
       No CL learning rate  & (1e-3, 2000, 0.95) & (1e-4, 2000, 0.99) & -- &  (1e-3, 2000, 0.99) \\
       No CL activation function & \texttt{swish} & \texttt{tanh} & -- & \texttt{relu}\\
       No CL architecture & [1, 200, 200, 200, 1] & [2, 200, 200, 200, 200, 200, 1] & -- & [5, 40, 40, 40, 1] \\
       No CL batch size  & 200 & 500 & -- & 100 \\
       No CL boundary condition batch size  & 1 & 100 & -- & -- \\
       No CL iterations  & 50,000 & 100,000 & -- & 100,000 \\
       $\lambda_{bc}$ & -- & 1 & -- & -- \\
       $\lambda_{ic}$ & 1 & 100 & -- & -- \\
       $\lambda_{r}$ & 10 & 1 & -- & -- \\
       $\lambda_{data}$ & 0 & 0 & 1 & 1 \\
       $\lambda$ & 1e-4 & 1e-5 & 1e-4 & 1e-5 \\
    \end{tabular}
    \caption{Hyperparameters for training the results in this paper. SF refers to the single fidelity results. MF refers to the multifidelity results. For the multifidelity results, $\mathcal{NN}^*$ denotes the first single fidelity network in the multifidelity framework, and $\mathcal{NN}_i$ denotes the multifidelity neural networks. For the learning rate, the triplet $(a, b, c)$ denotes the \texttt{exponential\_decay} function in Jax with learning rate $a$, decay steps $b$, and decay rate $c$. \texttt{relu} is the rectified linear unit (ReLU) activation function}
    \label{tab:parametesr}
\end{table}
\end{landscape}

\end{document}